\documentclass[12pt,a4paper,reqno]{amsart}

\usepackage{amsfonts}
\usepackage{amsmath}
\usepackage{amssymb}
\usepackage{amsthm}
\usepackage{amsxtra}
\usepackage{bbm}
\usepackage{braket}
\usepackage{comment}
\usepackage{extarrows}
\usepackage{graphicx}
\usepackage{latexsym}
\usepackage{mathrsfs}
\usepackage{yhmath}
\usepackage{nicematrix}
\usepackage{mathtools}

\usepackage{float}
\usepackage{booktabs}
\usepackage{verbatim}
\usepackage{xcolor}

\usepackage[pdfborder={0 0 0}]{hyperref} 
\hypersetup{breaklinks, raiselinks}

\usepackage{bm}

\usepackage{cleveref}

\usepackage{a4wide}
\usepackage{fullpage}

\theoremstyle{plain}
\newtheorem{Theorem}{Theorem}[section]
\newtheorem{Lemma}[Theorem]{Lemma}
\newtheorem{Corollary}[Theorem]{Corollary}
\newtheorem{Proposition}[Theorem]{Proposition}

\theoremstyle{definition}

\newtheorem{Assumptions and Discussion}[Theorem]{Assumptions and Discussion}

\newtheorem{Problem}[Theorem]{Problem}

\theoremstyle{remark}

\newtheorem*{acknowledgment*}{Acknowledgment}

\numberwithin{equation}{section}

\usepackage[shortlabels]{enumitem}
\SetEnumerateShortLabel{a}{\textup{(\alph*)}}
\SetEnumerateShortLabel{A}{\textup{(\Alph*)}}
\SetEnumerateShortLabel{1}{\textup{(\arabic*)}}
\SetEnumerateShortLabel{i}{\textup{(\roman*)}}
\SetEnumerateShortLabel{I}{\textup{(\Roman*)}}

\def\bar#1{\overline{#1}}

\begin{document}

\title{On Chern minimal surfaces in Hermitian surfaces}

\author[Chiakuei Peng, Xiaowei Xu]{Chiakuei Peng \quad  Xiaowei Xu${}^*$}





\address{School of Mathematical Sciences, University of Chinese Academy of Sciences, Beijing, 100049, P.R.~China.}
\email{pengck@ucas.edu.cn}

\address{CAS Wu Wen-Tsun Key Laboratory of Mathematics, School of Mathematical Sciences, University of Science and Technology of China, Hefei, Anhui, 230026, P.R.~China.}
\email{xwxu09@ustc.edu.cn}

\begin{abstract}
In this paper we introduce the Chern minimal surface in Hermitian surfaces by using the Chern connection, and we show that it only has isolated complex and anticomplex points for a generic one (neither holomorphic nor antiholomorphic). For a generic Chern minimal $f$ from compact Riemann surface $\Sigma$ in a Hermitian surface $M$, we establish two identities which related to the sum of the orders of all complex points, anticomplex points denoted by $P$, $Q$ respectively, the cap product of the pull-back of the first Chern class $f^*c_1(M)$ and $[\Sigma]$, the Euler characteristic of tangent bundle $\chi(T\Sigma)$ and the Euler characteristic of normal bundle $\chi(T^\perp\Sigma)$. More precisely, we obtain the formulae  $P-Q=-f^*c_1(M)[\Sigma]$ and $P+Q=-\big(\chi(T\Sigma)+\chi(T^\perp\Sigma)\big)$. We also give some applications of these formulae.
\end{abstract}
\maketitle

\section{Introduction}

The investigation on geometry and topology of minimal surfaces in four-dimensional manifold gained many geometer's interest. R.Bryant \cite{[B-82]} established the corresponding between superminimal surfaces into unit 4-sphere with holomorphic horizontal curves from surfaces into 3-dimensional complex projective space; S.S.Chern and J.G.Wolfson \cite{[CW-83]} found an important invariant, called K$\ddot{a}$hler angle, associated to surfaces immersed in a K$\ddot{a}$hlerian manifold; S.Webster \cite{[W-84]}, \cite{[W-86]} obtained two formulae with respect to the numbers of complex points, anticomplex points and topology data of minimal surfaces in a K$\ddot{a}$hler surface; J.H.Eschenburg, I.V.Guadalupe and R.A.Tribuzy \cite{[EGT-85]} studied the minimal surfaces in complex projective plane through the fundamental equations; J.G.Wolfson \cite{[W-89]} utilized the K$\ddot{a}$hler angle to give a new proof and global applications of the Webster formulae; M.J.Micallef and J.G.Wolfson \cite{[MW-93]} estimated the sum of index and nullity of Jacobi operator of minimal surfaces in a Riemannian 4-manifold by using the restriction on Weyl curvature, genus and normal Euler number of the surfaces; J.Chen and  G.Tian \cite{[CT-97]} defined an invariant, called the adjunction number to study the admissible surfaces in 4-dimensional Riemannian surfaces, and they also generalized one of Webster formulae to minimal surfaces in symplectic 4-manifold. In the paper presented, we wish to establish the Webster type formulae to Chern minimal surfaces in Hermitian surfaces.

Let $f$ be an immersion from Riemann surface $(\Sigma, \texttt{j})$ into K$\ddot{a}$hlerian or Hermitian surface $(M,J,g)$, where $\texttt{j}$, $J$ are the complex structures on $\Sigma$, $M$, respectively. Let $x$ be a point in $\Sigma$, we call $x$ is a complex point if $df\circ\texttt{j}=J\circ df$ holds at $x$, and we call $x$ is an anticomplex point if $df\circ\texttt{j}=-J\circ df$ holds at $x$, otherwise we call $x$ is a real point. Alternatively, by using the K$\ddot{a}$hler angle $\alpha:\Sigma\longrightarrow [0,\pi]$ defined by S.S.Chern and J.G.Wolfson in \cite{[CW-83]}, we know that $x\in\Sigma$ is a complex point (resp.anticomplex point) if and only if $\alpha(x)=0$ (resp. $\alpha(x)=\pi$), and $x$ is a real point if and only if $\alpha(x)\in (0,\pi)$. S.S.Chern and J.G.Wolfson show that the complex and anticomplex points are isolated if $f$ is  generic minimal (neither holomorphic nor antiholomorphic), and they also show that $\alpha$ is continuous on $\Sigma$, smooth away from the isolated complex and anticomplex point. Let $P$ be the number of complex points and $\;Q\;$ be the anticomplex points counted with 
the multiplicity. For a compact

\noindent------------------------------------------

*The corresponding author.

\noindent  Riemann surface immersed minimally into a K$\ddot{a}$hler surface $(M,J,g)$, 
 S.Webster \cite{[W-84]} proved that
\begin{equation}\label{1.1}
P+Q=-\big(\chi(T\Sigma)+\chi(T^\perp\Sigma)\big),
\end{equation}
where $\chi(T\Sigma)$ is the Euler characteristic of tangent bundle $T\Sigma$ and $\chi(T^\perp\Sigma)$ is the Euler characteristic of normal bundle $T^\perp\Sigma$; S.Webster \cite{[W-86]} also proved that
\begin{equation}\label{1.2}
P-Q=-f^*c_1(M)[\Sigma],
\end{equation}
where $c_1(M)$ is the first Chern class of $M$ and $[\Sigma]$ is the fundamental class of $\Sigma$. Now, the identities (\ref{1.1}) and (\ref{1.2}) are called Webster formulae. J.G.Wolfson \cite{[W-89]} give another proof of Webster formulae by considering the functions $\log\tan(\alpha/2)$ and $\log\sin\alpha$. Since the numbers $P$ and $Q$ on the left hand side of Webster formulae are only depend on the complex structures and tangent map, it is natural to ask the following
\begin{Problem}
\emph{Do the Webster formulae hold for certain class of compact Riemann surfaces immersed in Hermitian surfaces}?
\end{Problem}

Notice that the Chern connection preserves the metric, complex structure and has vanishing $(1,1)$-part of its torsion, and also the identity (\ref{1.2}) involving the first Chern class, 
it seems to find a condition related to the Chern connection can give an answer to this problem. For these reason, we first define the Chern minimal immersion from Riemann surface into Hermitian surfaces, which is trace-free of the second fundamental form induced by Chern connection. The Chern minimal is just minimal when the Hermitian surface is K$\ddot{a}$hlerian (see the Proposition 2.2), and the holomorphic and antiholomorphic immersions are automatically Chern minimal. For the generic Chern minimal surface in Hermitian surfaces, we prove

\begin{Theorem}\label{Theorem1.2}
Let $f$ be a generic Chern minimal immersion from the compact Riemann surface $\Sigma$ into Hermitian surface $M$. Then the Webster formulae $(\ref{1.1})$ and $(\ref{1.2})$ hold.
\end{Theorem}
The proof of Theorem 1.2 is much more sophisticated than that of K$\ddot{a}$hler surfaces, since the all calculations always involving the terms of torsion, however, they are vanishing for K$\ddot{a}$hler surfaces. For example, we need to construct the $1$-form $\theta_L$ which depend on the torsion of Chern connection, and  using the pull-back $f^* d\theta_L$ to measure the difference of $\Delta\log\tan(\alpha/2)\; dA$ and $f^*\mbox{Ric}_M$. This phenomenon
does not appeared in the K$\ddot{a}$hlerian case. Meanwhile, by reading the proof one can find that the Chern minimal condition is essential.

There are geometrical and topological applications of  Theorem 1.2 as following:

\begin{Theorem}\label{Theorem1.3}
Let $f$ be a generic Chern minimal immersion from the compact Riemann surface $\Sigma$ into Hermitian surface $M$. Then $f$ is real with constant K$\ddot{a}$hler angle if and only if $K+K^\perp=0$.
\end{Theorem}

\begin{Theorem}\label{Theorem1.4}
Let $f$ be a generic Chern minimal immersion from the compact Riemann surface $\Sigma$ into compact Hermitian surface $M$. Then, we have
\begin{equation*}
(2-2g)+\big|f^*c_1(M)[\Sigma]\big|+I_f-2D_f\leq -2\min\{P,Q\}\leq 0,
\end{equation*}
where $g$ is the genus of $\Sigma$, $I_f:=\langle D_M^{-1}f_*[\Sigma]\smile D_M^{-1}f_*[\Sigma], [M]\rangle$ is the intersection number of $\Sigma$, $D_f$ is the self-intersection number which  has only regular self intersections of multiplicity $2$.
\end{Theorem}
\emph{Remak}. We should to point out that the K$\ddot{a}$hlerian version of Theorem 1.4 is due to J.G.Wolfson \cite{[W-89]}. There are some interesting corollaries of Theorem 1.4, one can refer to the Section 5 for details.

The paper organized as follows. In section 2, we set up the basic formulae for smooth immersions from Riemann surface into Hermitian surfaces, we present the definition of Chern minimal and show the difference between Chern minimal and usual one, and we also recall the definition of K$\ddot{a}$hler angle; In Section 3, 4, we prove the identities (\ref{1.2}), (\ref{1.1}) respectively; In Section 5, we give some applications of our main theorem.

\vspace{0.6cm}

\section{Chern Minimal and K$\ddot{\mathrm{A}}$hler Angle}
The purpose of this section is to establish notations, the basic formulae of Riemann surfaces isometrically immersed in a Hermitian surfaces related to the Chern connection.

\subsection{Surfaces in Hermitian Surface}
Let $(\Sigma,\texttt{j})$ be a Riemann surface, and let $ds^2_\Sigma$ be a Riemannian metric on $\Sigma$ which is conformal to the complex structure 
$\texttt{j}$. Locally, there exists a complex valued 1-form $\varphi$ of $(1,0)$-type such that 
\begin{equation}\label{2.1}
ds^2_\Sigma=\varphi\;\bar{\varphi}.
\end{equation}
The first structure
equation of $ds^2_\Sigma$ is given by
\begin{equation}\label{2.2}
d\varphi=-\sqrt{-1}\rho\wedge \varphi,
\end{equation}
where $\rho$ is the connection 1-form. The Gaussian curvature $K$ is determined by 
\begin{equation}\label{2.3}
d\rho=-\sqrt{-1}K \varphi\wedge \bar{\varphi}.
\end{equation}

Let $(M,J,g)$ be a Hermitian surface, where $J$ is the complex structure and $g$ is the Hermitian metric. Let $\nabla$ be the Chern connection on the tangent bundle $TM$, which is the unique connection with vanishing $(1,1)$-part of torsion and preserving $J$ and $g$, i.e., 
\begin{equation}\label{2.4}
\nabla J=0,\hspace{0.3cm} \nabla g=0.
\end{equation}
More explicitly, choosing a local unitary frame $\{e_1,e_2\}$ with the dual $\{\omega^1,\omega^2\}$. Then, the fundamental 2-form $\omega$ is given by
\begin{equation}\label{2.5}
\omega=\sqrt{-1}(\omega^1\wedge\overline{\omega^1}+\omega^2\wedge\overline{\omega^2}).
\end{equation}
 It follows from (\ref{2.4}) that the \emph{connection $1$-forms} $\{\omega^j_i\}$ are given by
\begin{equation}\label{2.6}
\nabla e_i=\omega_i^j e_j,\hspace{1cm}\omega^j_i+\overline{\omega^i_j}=0.
\end{equation}
The first structure equations of $\nabla$ are determined by
\begin{equation}\label{2.7}
d\omega^i=-\omega^i_j\wedge\omega^j+\Theta^i,
\end{equation}
where the torsion $\Theta^i$ can be written as
\begin{equation}\label{2.8}
\Theta^i:=L^i_{jk}\,\omega^j\wedge\omega^k, \hspace{0.3cm}L^i_{jk}=-L^i_{kj}.
\end{equation}
The second structure equations of $\nabla$ are determined by 
\begin{equation}\label{2.9}
d\omega^i_j=-\omega^i_k\wedge\omega^k_j+\Omega^i_j,\hspace{0.5cm}\Omega^j_i+\overline{\Omega^i_j}=0,
\end{equation}
where the curvature $2$-forms $\Omega^i_j$ can be written as
\begin{equation}\label{2.10}
\Omega^i_j:=R^i_{jk\ell}\,\omega^k\wedge\omega^\ell+R^i_{jk\bar{\ell}}\,\omega^k\wedge\overline{\omega^\ell}+R^i_{j\bar{k}\,\bar{\ell}}\,\overline{\omega^k}\wedge\overline{\omega^\ell}
\end{equation}
with $R^i_{jk\ell}=-R^i_{j\ell k}$, $R^i_{j\bar{k}\,\bar{\ell}}=-R^i_{j\bar{\ell}\,\bar{k}}$. The skew-Hermitian of curvature forms $\Omega^i_j$ implies
\begin{equation}\label{2.11}
R^i_{jk\ell}=\overline{R^j_{i\bar{\ell}\,\bar{k}}},\hspace{1cm}R^i_{jk\bar{\ell}}=\overline{R^j_{i\ell\bar{k}}}.
\end{equation}

Throughout this paper we always assume that $f:\Sigma\longrightarrow M$ is a smooth conformal isometric immersion from Riemann surface $\Sigma$ into Hermitian surface $(M,J,g)$. Set
\begin{equation}\label{2.12}
f^*\omega^i:=a^i_1\,\varphi+a^i_{\bar{1}}\,\bar{\varphi}.
\end{equation}
Then the conformality and isometricity of $f$ are equivalent to
\begin{equation}\label{2.13}
a^1_1\,\overline{a^1_{\bar{1}}}+a^2_1\,\overline{a^2_{\bar{1}}}=0,
\end{equation}
\begin{equation}\label{2.14}
|a^1_1|^2+|a^2_1|^2+|a^1_{\bar{1}}|^2+|a^2_{\bar{1}}|^2=1,
\end{equation}
respectively. Taking the exterior derivative of (\ref{2.12}) and using the structure equations (\ref{2.2}) and (\ref{2.7}), we have
\begin{equation}\label{2.15}
(a^i_{11}\,\varphi+a^i_{1\bar{1}}\,\bar{\varphi})\wedge\varphi+(a^i_{\bar{1}1}\,\varphi+a^i_{\bar{1}\,\bar{1}}\,\bar{\varphi})\wedge\bar{\varphi}=f^*\Theta^i,
\end{equation}
where $a^i_{11}$, $a^i_{1\bar{1}}$, $a^i_{\bar{1}1}$, $a^i_{\bar{1}\,\bar{1}}$ are defined as 
\begin{equation}\label{2.16}
a^i_{11}\,\varphi+a^i_{1\bar{1}}\,\bar{\varphi}:=da^i_1-\sqrt{-1}\rho\, a^i_1 +a^j_1\,\omega^i_j,
\end{equation}
\begin{equation}\label{2.17}
a^i_{\bar{1}1}\,\varphi+a^i_{\bar{1}\,\bar{1}}\,\bar{\varphi}:=da^i_{\bar{1}}+\sqrt{-1}\rho\, a^i_{\bar{1}} +a^j_{\bar{1}}\,\omega^i_j.
\end{equation}
Here, we have omitted the pull-back $f^*$ acts on the forms defined on $M$. By comparing the coefficients of (\ref{2.15}) and using (\ref{2.8}), (\ref{2.12}), we have
\begin{equation}\label{2.18}
-a^i_{1\bar{1}}+a^i_{\bar{1}1}=2L^i_{jk}a^j_1a^k_{\bar{1}}.
\end{equation}
Similarly, taking the exterior derivative of (\ref{2.16}), (\ref{2.17}), we obtain
\begin{equation}\label{2.19}
(a^i_{111}\varphi+a^i_{11\bar{1}}\,\bar{\varphi})\wedge\varphi+(a^i_{1\bar{1}1}\varphi+a^i_{1\bar{1}\,\bar{1}}\,\bar{\varphi})\wedge\bar{\varphi}=-Ka^i_1\varphi\wedge\bar{\varphi}+a^j_1\Omega_j^i,
\end{equation}
\begin{equation}\label{2.20}
(a^i_{\bar{1}11}\varphi+a^i_{\bar{1}1\bar{1}}\,\bar{\varphi})\wedge\varphi+(a^i_{\bar{1}\,\bar{1}1}\varphi+a^i_{\bar{1}\,\bar{1}\,\bar{1}}\,\bar{\varphi})\wedge\bar{\varphi}=Ka^i_{\bar{1}}\,\varphi\wedge\bar{\varphi}+a^j_{\bar{1}}\,\Omega_j^i,\hspace{0.1cm}{}
\end{equation}
respectively, where $a^i_{111}$, $a^i_{11\bar{1}}$, $a^i_{1\bar{1}1}$, $a^i_{1\bar{1}\,\bar{1}}$, $a^i_{\bar{1}11}$, $a^i_{\bar{1}1\bar{1}}$, $a^i_{\bar{1}\,\bar{1}1}$, $a^i_{\bar{1}\,\bar{1}\,\bar{1}}$ are defined by
\begin{equation}\label{2.21}
a^i_{111}\varphi+a^i_{11\bar{1}}\,\bar{\varphi}:=da^i_{11}-2a^i_{11}\sqrt{-1}\rho+a^j_{11}\omega^i_j,
\end{equation}
\begin{equation}\label{2.22}
a^i_{1\bar{1}1}\varphi+a^i_{1\bar{1}\,\bar{1}}\,\bar{\varphi}:=da^i_{1\bar{1}}+a^j_{1\bar{1}}\,\omega^i_j,
\end{equation}
\begin{equation}\label{2.23}
a^i_{\bar{1}11}\varphi+a^i_{\bar{1}1\bar{1}}\,\bar{\varphi}:=da^i_{\bar{1}1}+a^j_{\bar{1}1}\omega^i_j,
\end{equation}
\begin{equation}\label{2.24}
a^i_{\bar{1}\,\bar{1}1}\varphi+a^i_{\bar{1}\,\bar{1}\,\bar{1}}\,\bar{\varphi}:=da^i_{\bar{1}\,\bar{1}}+2a^i_{\bar{1}\,\bar{1}}\sqrt{-1}\rho+a^j_{\bar{1}\,\bar{1}}\omega^i_j.
\end{equation}
By comparing the coefficients of (\ref{2.19}), (\ref{2.20}), we obtain
\begin{equation}\label{2.25}
a^i_{1\bar{1}1}-a^i_{11\bar{1}}=-Ka^i_1+a^j_1\,\Omega^i_j/(\varphi\wedge\bar{\varphi}),
\end{equation}
\begin{equation}\label{2.26}
a^i_{\bar{1}\,\bar{1}1}-a^i_{\bar{1}1\bar{1}}=Ka^i_{\bar{1}}+a^j_{\bar{1}}\,\Omega^i_j/(\varphi\wedge\bar{\varphi}),
\end{equation}
where $\Omega^i_j/(\varphi\wedge\bar{\varphi})$ stands for the coefficient of the pull-back of $\Omega^i_j$ with respect to $\varphi\wedge\bar{\varphi}$.

It follows from (\ref{2.12}) that the differential of $f$ can be expressed by
\begin{equation*}
df=a^i_1\,\varphi\otimes e_i+a^i_{\bar{1}}\,\bar{\varphi}\otimes e_i
+\overline{a^i_1}\,\overline{\varphi}\otimes\overline{e_i}+\overline{a^i_{\bar{1}}}\,\varphi\otimes\overline{e_i},
\end{equation*}
which is a smooth section of bundle $T^*\Sigma\otimes f^{-1}TM$ or $(T^{\mathbb{C}}\Sigma)^*\otimes f^{-1}T^{\mathbb{C}}M$. Taking the covariant derivative of $df$ by using the induced connection, we get the \emph{second fundamental form} $\nabla df$. Explicitly, by using (\ref{2.16}), (\ref{2.17}), we have
\begin{eqnarray*}
\nabla df&=&a^i_{11}\,\varphi\otimes\varphi\otimes e_i+a^i_{1\bar{1}}\,\bar{\varphi}\otimes\varphi\otimes e_i
+a^i_{\bar{1}1}\,\varphi\otimes\bar{\varphi}\otimes e _i+a^i_{\bar{1}\,\bar{1}}\,\bar{\varphi}\otimes\bar{\varphi}\otimes e_i\nonumber\\
&{}&+\,\overline{a^i_{11}}\,\bar{\varphi}\otimes\bar{\varphi}\otimes \overline{e_i}+\overline{a^i_{1\bar{1}}}\,\varphi\otimes\bar{\varphi}\otimes\overline{e_i}
+\overline{a^i_{\bar{1}1}}\,\bar{\varphi}\otimes\varphi\otimes \overline{e_i}+\overline{a^i_{\bar{1}\,\bar{1}}}\,\varphi\otimes\varphi\otimes \overline{e_i}.
\end{eqnarray*}
Since $tr\nabla df$ is the mean curvature of $f$ which is closely related to the Chern connection, we call it is \emph{Chern minimal} if $tr\nabla df=0$, and it is equivalent to 
\begin{equation}\label{2.27}
a^i_{1\bar{1}}+a^i_{\bar{1}1}=0,
\end{equation}
for all $i$. It is easy to check that a holomorphic or antiholomorphic isometric immersion is Chern minimal. So, it is natural to call a Chern minimal immersion $f$ is \emph{generic Chern minimal} if $f$ is neither holomorphic nor antiholomorphic.

\subsection{Levi-Civita Connection}
Denote the Levi-Civita connection of Hermitian surface $(M,J,g)$ by $D$. We extend the Levi-Civita connection $\mathbb{C}$-linearly to the complexified tangent bundle $T^\mathbb{C}M$, and set
\begin{equation}\label{2.28}
De_i:=\varphi_i^j\,e_j+\varphi^{\bar{j}}_i\,\overline{e_j},
\end{equation}
where $\varphi_i^j$, $\varphi^{\bar{j}}_i$ are complex-valued 1-forms. It follows from $Dg=0$ that we have
\begin{equation}\label{2.29}
\varphi_i^j+\overline{\varphi^i_j}=0,\hspace{1cm}\varphi_i^{\bar{j}}+\varphi_j^{\bar{i}}=0.
\end{equation}
It is easy to check that the first structure equation of $D$ is given by
\begin{equation}\label{2.30}
d\omega^i=-\varphi^i_j\wedge\omega^j-\overline{\varphi_j^{\bar{i}}}\wedge\overline{\omega^j}.
\end{equation}

\begin{Proposition}
The differences between the connections of $\nabla$ and $D$ are given by 

\emph{(1)}  $\varphi^i_j=\omega^i_j+L^i_{jk}\,\omega^k-\overline{L^j_{ik}}\,\overline{\omega^k}.$

\emph{(2)}  $\varphi^{\bar{j}}_i=L^k_{ij}\,\overline{\omega^k}$.
\end{Proposition}

\emph{Proof}. It follows the first structure equations (\ref{2.7}) and (\ref{2.30}) we obtain
\begin{equation}\label{2.31}
(\varphi^i_j-\omega^i_j)\wedge\omega^j+\overline{\varphi^{\bar{i}}_j}\wedge\overline{\omega^j}+\Theta^i=0.
\end{equation}
For fixed indices $i,k,\ell$, evaluating (\ref{2.31}) on the pairs $(e_k,e_\ell)$, $(e_k,\overline{e_\ell})$, $(\overline{e_k},\overline{e_\ell})$ repectively, we have
\begin{equation}\label{2.32}
\big(\varphi^i_\ell(e_k)-\omega^i_\ell(e_k)\big)-\big(\varphi^i_k(e_\ell)-\omega^i_k(e_\ell)\big)+2L^i_{k\ell}=0,
\end{equation}
\begin{equation}\label{2.33}
\big(\varphi^i_k(\overline{e_\ell})-\omega^i_k(\overline{e_\ell})\big)-\overline{\varphi^{\bar{i}}_\ell}(e_k)=0,
\end{equation}
\begin{equation}\label{2.34}
\overline{\varphi^{\bar{i}}_\ell}(\overline{e_k})-\overline{\varphi^{\bar{i}}_k}(\overline{e_\ell})=0.
\end{equation}
Taking the conjugate of (\ref{2.33}) and using the skew-Hermitian of $\varphi^i_j$, $\omega^i_j$, we obtain
\begin{equation}\label{2.35}
\big(\varphi^k_i(e_\ell)-\omega^k_i(e_\ell)\big)+\varphi^{\bar{i}}_\ell(\overline{e_k})=0,
\end{equation}
which implies
\begin{equation}\label{2.36}
\varphi^{\bar{\ell}}_k(\overline{e_i})=L^i_{k\ell}
\end{equation}
by (\ref{2.32}) and the skew-symmetric of $\varphi^{\bar{\ell}}_k$. Then the first identity is follows from substituting (\ref{2.36}) into (\ref{2.33}) and (\ref{2.35}). On the other hand, circling the indices of the conjugation of (\ref{2.34}) to obtain three identities, summing any two of them and minus the other one, we can get
\begin{equation}\label{2.37}
\varphi^{\bar{\ell}}_k(e_i)=0.
\end{equation}
So, the second identity follows from (\ref{2.36}) and (\ref{2.37}). 

\hfill$\Box$

\vspace{0.3cm}
\begin{Proposition}\label{Proposition 2.2}
The differences between of $trDdf$ and $tr\nabla df$ are given by
\begin{eqnarray*}
tr Ddf=tr\nabla df+2(a^j_{\bar{1}}\,\overline{L^j_{ki}}\,\overline{a^k_{\bar{1}}}+\overline{a^j_1}\,\overline{L^k_{ji}}\,a^k_1)\,e_i
+2(\overline{a^j_{\bar{1}}}\,L^j_{ki}\,a^k_{\bar{1}}+a^j_1\,L^k_{ji}\,\overline{a^k_1})\,\overline{e_i}.
\end{eqnarray*}
\end{Proposition}

\emph{Proof}. By using the differences between $D$ and $\nabla$ given in the Proposition 2.1, we have
\begin{eqnarray} \label{2.38}
D(a^i_1\,\varphi\otimes e_i)&=&(a^i_{11}-a^j_1\,\overline{L^j_{ik}}\,\overline{a^k_{\bar{1}}})\;\varphi\otimes\varphi\otimes e_i
+(a^i_{1\bar{1}}+a^j_1\,L^i_{jk}\,a^k_{\bar{1}}-a^j_1\,\overline{L^j_{ik}}\,\overline{a^k_1})\;\bar{\varphi}\otimes\varphi\otimes e_i\nonumber\\
&{}&+a^j_1\,L^k_{ji}\,\overline{a^k_{\bar{1}}}\;\varphi\otimes\varphi\otimes\overline{e_i}+a^j_1\,L^k_{ji}\, \overline{a^k_1}\;\bar{\varphi}\otimes\varphi\otimes\overline{e_i},
\end{eqnarray}
and 
\begin{eqnarray} \label{2.39}
D(a^i_{\bar{1}}\,\bar{\varphi}\otimes e_i)&=&(a^i_{\bar{1}1}+a^j_{\bar{1}}\,L^i_{jk}\,a^k_1-a^j_{\bar{1}}\,\overline{L^j_{ik}}\,\overline{a^k_{\bar{1}}})\;\varphi\otimes\bar{\varphi}\otimes e_i
+(a^i_{\bar{1}\,\bar{1}}-a^j_{\bar{1}}\,\overline{L^j_{ik}}\,\overline{a^k_1})\;\bar{\varphi}\otimes\bar{\varphi}\otimes e_i\nonumber\\
&{}&+a^j_{\bar{1}}\,L^k_{ji}\,\overline{a^k_{\bar{1}}}\;\varphi\otimes\bar{\varphi}\otimes\overline{e_i}+a^j_{\bar{1}}\,L^k_{ji}\, \overline{a^k_1}\;\bar{\varphi}\otimes\bar{\varphi}\otimes\overline{e_i}.
\end{eqnarray}
Substuting (\ref{2.38}), (\ref{2.39}) and their conjugate into
\begin{equation*}
Ddf=D(a^i_1\,\varphi\otimes e_i+a^i_{\bar{1}}\,\bar{\varphi}\otimes e_i
+\overline{a^i_1}\,\overline{\varphi}\otimes\overline{e_i}+\overline{a^i_{\bar{1}}}\,\varphi\otimes\overline{e_i}),
\end{equation*}
we obtain
\begin{eqnarray}\label{2.40}
Ddf\!\!&=&\!\!\nabla df-2a^j_1\,L^j_{ik}\,\overline{a^k_{\bar{1}}}\,\varphi\otimes\varphi\otimes e_i+(a^j_{\bar{1}}\,L^i_{jk}\,a^k_1-a^j_{\bar{1}}\,\overline{L^j_{ik}}\,\overline{a^k_{\bar{1}}}+\overline{a^j_1}\,\overline{L^k_{ji}}\,a^k_1)\varphi\otimes\bar{\varphi}\otimes e_i
\nonumber\\
&{}&+(a^j_1\,L^i_{jk}\,a^k_{\bar{1}}-a^j_1\,\overline{L^j_{ik}}\,\overline{a^k_1}+\overline{a^j_{\bar{1}}}\,\overline{L^k_{ji}}\,a^k_{\bar{1}})\,\bar{\varphi}\otimes\varphi\otimes e_i-2a^j_{\bar{1}}\,\overline{L^j_{ik}}\,\overline{a^k_1}\,\bar{\varphi}\otimes\bar{\varphi}\otimes e_i
\nonumber\\
&{}&+2a^j_1\,L^k_{ji}\,\overline{a^k_{\bar{1}}}\,\varphi\otimes\varphi\otimes \overline{e_i}+(a^j_{\bar{1}}\,L^k_{ji}\,\overline{a^k_{\bar{1}}}+\overline{a^j_1}\,\overline{L^i_{jk}}\,\overline{a^k_{\bar{1}}}-\overline{a^j_1}\,L^j_{ik}\,a^k_1)\,\varphi\otimes\bar{\varphi}\otimes \overline{e_i}\nonumber\\
&{}&+(a^j_1\,L^k_{ji}\,\bar{a^k_1}+\overline{a^j_{\bar{1}}}\,\overline{L^i_{jk}}\,\overline{a^k_1}-\overline{a^j_{\bar{1}}}\,L^j_{ik}\,a^k_{\bar{1}})\bar{\varphi}\otimes\varphi\otimes\overline{e_i}+2a^j_{\bar{1}}\,L^k_{ji}\,\overline{a^k_1}\,\bar{\varphi}\otimes\bar{\varphi}\otimes\overline{e_i}.
\end{eqnarray}
Denote the dual of $\varphi$ by $e$, then (\ref{2.40}) gives 
\begin{eqnarray}\label{2.41}
trDdf&=&Ddf(e,\bar{e})+Ddf(\bar{e},e)\nonumber\\
&=& tr\nabla df+2(a^j_{\bar{1}}\,\overline{L^j_{ki}}\,\overline{a^k_{\bar{1}}}+\overline{a^j_1}\,\overline{L^k_{ji}}\,a^k_1)\,e_i
+2(\overline{a^j_{\bar{1}}}\,L^j_{ki}\,a^k_{\bar{1}}+a^j_1\,L^k_{ji}\,\overline{a^k_1})\,\overline{e_i}.
\end{eqnarray}
We have used the fact $L^i_{jk}=-L^i_{kj}$ to get (\ref{2.40}) and (\ref{2.41}).

\hfill$\Box$

\emph{Remark}. It is known that $f$ is minimal if and only if $trDdf=0$. So, this proposition tells us that the ``minimal'' and ``Chern minimal'' are different in general. But, they are coincide with each other when $(M,J,g)$ is K$\ddot{a}$hlerian.
We will focus on the Chern minimal surfaces in the sequal.

\subsection{K$\ddot{a}$hler Angle}
As S.S.Chern and J.G.Wolfson pointed out in \cite{[CW-83]}, there is an important invariant associated to a conformal isometric immersion $f$ from Riemann surface $\Sigma$ into Hermitian surface $M$, called K$\ddot{a}$hler angle, which measures the holomorphicity of $f$. Locally, by the conformality condition (\ref{2.13}), we can choose a unitary coframe $\{\omega^1,\omega^2\}$ for $M$ such that 
\begin{equation}\label{3.1}
\omega^1=s\,\varphi,\hspace{0.3cm}\omega^2=t\,\bar{\varphi},
\end{equation}
where $s$, $t$ are local smooth complex valued fuctions on $\Sigma$. By the isometricity condition, we can set
\begin{equation}\label{3.2}
|s|=\cos \frac{\alpha}{2},\hspace{0.3cm} |t|=\sin \frac{\alpha}{2},
\end{equation}
where $\alpha$ is a function on $\Sigma$ with values between $0$ and $\pi$, which is called \emph{K$\ddot{a}$hler angle}. At a point $x\in\Sigma$, the tangent space $T_x\Sigma$ is a \emph{complex} (resp. \emph{anticomplex}) line in $T_{f(x)}M$ when $\alpha(x)=0$ (resp. $\alpha(x)=\pi$), and such points are called \emph{complex} (resp. \emph{anticomplex}) \emph{points}. Otherwise, we call $x\in\Sigma$ is a \emph{real} point if $0<\alpha(x)<\pi$, and call it is a \emph{totally real} point if $\alpha(x)=\pi/2$ in particular. We call $f$ is \emph{holomorphic} (resp. \emph{antiholomorphic}) if every point in $\Sigma$ is complex (resp. anticomplex), and similar for $f$ is \emph{real} and \emph{totally real}. We should remark that $\alpha$ is a smooth function away from complex and anticomplex points, and it is continuous at these points.

\begin{Proposition} \label{Proposition 2.3}
Let $f$ be a generic Chern minimal immersion from Riemann surface $\Sigma$ into Hermitian surface $M$. Then the complex  points and anticomplex points of $f$ are isolated.
\end{Proposition}

\emph{Proof}. It follows form (\ref{3.1}) we have $a^1_1=s$ and $a^2_1=0$. Substuting them into (\ref{2.16}) we obtain
\begin{equation*}
ds=s(\sqrt{-1}\rho-\omega^1_1)+a^1_{11}\varphi+a^1_{1\bar{1}}\bar{\varphi}.
\end{equation*}
On the other hand, the Chern minimality (\ref{2.27}) and (\ref{2.18}), we obtain $a^1_{1\bar{1}}=-stL^1_{12}$. Thus, we have
\begin{equation}\label{3.3}
ds=s(\sqrt{-1}\rho-\omega^1_1-tL^1_{12}\bar{\varphi})+a^1_{11}\varphi.
\end{equation}
So, in a holomorphic coordinate $z$ around a point $x\in \Sigma$, (\ref{3.3}) implies $s$ satisfies the equation
\begin{equation*}
\frac{\partial s}{\partial \bar{z}}=s\,u,
\end{equation*}
where $u$ is the coefficient of $d\bar{z}$ in $\sqrt{-1}\rho-\omega^1_1-tL^1_{12}\bar{\varphi}$. By the main theorem in Section 4 of \cite{[Ch-70]}, $s$ can be written as
\begin{equation}\label{3.4}
s=z^p\,h,
\end{equation}
where $p\in\textbf{N}$ and $h$ is a function without zero point. Thus, $s$ only has isolated zeroes, and similar proof for $t$.

\hfill$\Box$

\emph{Remark}. For an isolated zero $x$ of $s$, the positive integer $p$ in (\ref{3.4}) is called the \emph{order} of $x$, which is independent of the choice of holomorphic coordinate.

Let $f$ be a generic Chern minimal immersion from compact Riemann surface $\Sigma$ into Hermitian surface $M$, the complex points and antiholomorphic points are denoted by $\{x_i\,|\,i=1,\ldots,m\}$ and $\{y_j\,|\,j=1,\dots,n\}$, and the orders of $x_i$, $y_j$ are denoted by $p_i$, $q_j$, respectively. Then all the complex points $P$ and anticomplex points $Q$ counted with orders are given by
\begin{equation}\label{3.5}
P=p_1+\ldots+p_m,\hspace{0.5cm}Q=q_1+\ldots+q_n,
\end{equation}
respectively.
By the local expression (\ref{3.4}) and Stokes' formula, it is easy to check  that
\begin{equation}\label{3.6}
P=-\frac{1}{2\pi}\int_\Sigma \Delta\log |t|\;dA,\hspace{0.5cm}Q=-\frac{1}{2\pi}\int_\Sigma \Delta\log |s|\;dA,
\end{equation}
where $\Delta$, $dA$ are the Laplacian, area element of $ds^2_\Sigma$ respectively.

\vspace{0.6cm}

\section{The Cap Product $f^*c_1(M)[\Sigma]$}

To begin, we do calculations where $\alpha$ is smooth. Locally, we can further choose a unitary frame $\{\omega^1, \omega^2\}$ satisfying
\begin{equation}\label{3.7}
\omega^1=\cos\frac{\alpha}{2}\,\varphi,\hspace{0.8cm}\omega^2=\sin\frac{\alpha}{2}\,\bar{\varphi},
\end{equation}
which imply
\begin{equation}\label{3.8}
a^1_1=\cos\frac{\alpha}{2},\hspace{0.5cm} a^1_1=0,\hspace{0.5cm} a^1_{\bar{1}}=0, \hspace{0.5cm}a^2_{\bar{1}}=\sin\frac{\alpha}{2}.
\end{equation}
Substituting them into (\ref{2.16}), (\ref{2.17}) and using the second identity in (\ref{2.6}), we obtain
\begin{equation}\label{3.9}
\frac{1}{2}\sin\frac{\alpha}{2}\,d\alpha-\sqrt{-1}\cos\frac{\alpha}{2}\,\rho+\cos\frac{\alpha}{2}\,\omega^1_1=-\overline{a^1_{1\bar{1}}}\,\varphi-\overline{a^1_{11}}\,\bar{\varphi},
\end{equation}
\begin{equation}\label{3.10}
\cos\frac{\alpha}{2}\,\omega^2_1=a^2_{11}\,\varphi+a^2_{1\bar{1}}\,\bar{\varphi},
\end{equation}
\begin{equation}\label{3.11}
\sin\frac{\alpha}{2}\,\omega^2_1=-\overline{a^1_{\bar{1}\,\bar{1}}}\,\varphi-\overline{a^1_{\bar{1}1}}\,\bar{\varphi},
\end{equation}
\begin{equation}\label{3.12}
\frac{1}{2}\cos\frac{\alpha}{2}\,d\alpha+\sqrt{-1}\sin\frac{\alpha}{2}\,\rho+\sin\frac{\alpha}{2}\,\omega^2_2=a^2_{\bar{1}1}\,\varphi+a^2_{\bar{1}\,\bar{1}}\,\bar{\varphi}.
\end{equation}
Then, the identities (\ref{3.9}), (\ref{3.12}) imply
\begin{equation}\label{3.13}
\frac{1}{2}\big[d\alpha+\sin\alpha\,(\omega^1_1+\omega^2_2)\big]=(-\sin\frac{\alpha}{2}\,\overline{a^1_{1\bar{1}}}+\cos\frac{\alpha}{2}\,a^2_{\bar{1}1})\,\varphi+(-\sin\frac{\alpha}{2}\,\overline{a^1_{11}}+\cos\frac{\alpha}{2}\,a^2_{\bar{1}\,\bar{1}})\,\bar{\varphi},
\end{equation}
and the identities (\ref{3.10}), (\ref{3.11}) imply
\begin{equation}\label{3.14}
\omega^2_1=(\cos\frac{\alpha}{2}\,a^2_{11}-\sin\frac{\alpha}{2}\,\overline{a^1_{\bar{1}\,\bar{1}}}\,)\,\varphi+
(\cos\frac{\alpha}{2}\,a^2_{1\bar{1}}-\sin\frac{\alpha}{2}\,\overline{a^1_{\bar{1}1}}\,)\,\bar{\varphi}.
\end{equation}
The identity (\ref{3.13}) and its conjugate give
\begin{eqnarray}\label{3.15}
d\alpha&=&(-\sin\frac{\alpha}{2}\,a^1_{11}-\sin\frac{\alpha}{2}\,\overline{a^1_{1\bar{1}}}+\cos\frac{\alpha}{2}\,a^2_{\bar{1}1}+\cos\frac{\alpha}{2}\,\overline{a^2_{\bar{1}\,\bar{1}}})\,\varphi\nonumber\\
&{}&-(\sin\frac{\alpha}{2}\,\overline{a^1_{11}}+\sin\frac{\alpha}{2}\,a^1_{1\bar{1}}-\cos\frac{\alpha}{2}\,\overline{a^2_{\bar{1}1}}-\cos\frac{\alpha}{2}\,a^2_{\bar{1}\,\bar{1}})\,\bar{\varphi},
\end{eqnarray}
and
\begin{eqnarray}\label{3.16}
\sin\alpha(\omega^1_1+\omega^2_2)&=&(\sin\frac{\alpha}{2}\,a^1_{11}-\sin\frac{\alpha}{2}\,\overline{a^1_{1\bar{1}}}+\cos\frac{\alpha}{2}\,a^2_{\bar{1}1}-\cos\frac{\alpha}{2}\,\overline{a^2_{\bar{1}\,\bar{1}}})\,\varphi\nonumber\\
&{}&-(\sin\frac{\alpha}{2}\,\overline{a^1_{11}}-\sin\frac{\alpha}{2}\,a^1_{1\bar{1}}+\cos\frac{\alpha}{2}\,\overline{a^2_{\bar{1}1}}-\cos\frac{\alpha}{2}\,a^2_{\bar{1}\,\bar{1}})\,\bar{\varphi}.
\end{eqnarray}
Substituting (\ref{3.14}) back to (\ref{3.10}) and (\ref{3.11}), by comparing the coefficients of $\varphi$, $\bar{\varphi}$ respectively, we have
\begin{equation}\label{3.17}
\sin^2\frac{\alpha}{2}\,a^2_{11}=-\cos\frac{\alpha}{2}\sin\frac{\alpha}{2}\,\overline{a^1_{\bar{1}\,\bar{1}}},
\end{equation}
\begin{equation}\label{3.18}
\sin^2\frac{\alpha}{2}\,a^2_{1\bar{1}}=-\cos\frac{\alpha}{2}\sin\frac{\alpha}{2}\,\overline{a^1_{\bar{1}\,1}},
\end{equation}
\begin{equation}\label{3.19}
\cos^2\frac{\alpha}{2}\,\overline{a^1_{\bar{1}\,\bar{1}}}=-\cos\frac{\alpha}{2}\sin\frac{\alpha}{2}\,a^2_{11},
\end{equation}
\begin{equation}\label{3.20}
\cos^2\frac{\alpha}{2}\,\overline{a^1_{\bar{1}1}}=-\cos\frac{\alpha}{2}\sin\frac{\alpha}{2}\,a^2_{1\,\bar{1}}.
\end{equation}

Notice that
$\cos^2\frac{\alpha}{2}=a^i_1\,\overline{a^i_1}$, $\sin^2\frac{\alpha}{2}=a^i_{\bar{1}}\,\overline{a^i_{\bar{1}}}$,
we have
\begin{eqnarray}\label{3.21}
\big(\cos^2\frac{\alpha}{2}\big)_{1}&=&\big(a^i_1\,\overline{a^i_1}\big)_1\nonumber\\
&=&a^i_{11}\,\overline{a^i_1}+a^i_1\,\overline{a^i_{1\bar{1}}}\nonumber\\
&=&\cos\frac{\alpha}{2}\,(a^1_{11}+\overline{a^1_{1\bar{1}}}),
\end{eqnarray}
and
\begin{eqnarray}\label{3.22}
\frac{1}{2}\Delta\cos^2\frac{\alpha}{2}&=&\big(a^i_1\,\overline{a^i_1}\big)_{1\bar{1}}\nonumber\\
&=&a^i_{11\bar{1}}\,\overline{a^i_1}+|a^i_{11}|^2+|a^i_{1\bar{1}}|^2+a^i_1\,\overline{a^i_{1\bar{1}1}}\nonumber\\
&=&(a^i_{1\bar{1}1}+Ka^i_1-a^j_1\,\Omega_j^i/\varphi\wedge\bar{\varphi})\,\overline{a^i_1}+|a^i_{11}|^2+|a^i_{1\bar{1}}|^2+a^i_1\,\overline{a^i_{1\bar{1}1}}\nonumber\\
&=&2\mbox{Re}(\overline{a^i_1}\,a^i_{1\bar{1}1})+|a^i_{11}|^2+|a^i_{1\bar{1}}|^2+\cos^2\frac{\alpha}{2}\,(K-\Omega^1_1/\varphi\wedge\bar{\varphi}),
\end{eqnarray}
where we have used the Ricci identity (\ref{2.25}). Substituting (\ref{3.21}) and its conjugate, (\ref{3.22}) into
\begin{eqnarray*}
\frac{1}{2}\Delta\log\cos^2\frac{\alpha}{2}&=&\big(\log\cos^2\frac{\alpha}{2}\big)_{1\bar{1}}\\
&=& -\frac{(\cos^2\frac{\alpha}{2})_1\;(\cos^2\frac{\alpha}{2})_{\bar{1}}}{\cos^4\frac{\alpha}{2}}+\frac{1}{\cos^2\frac{\alpha}{2}}\cdot\frac{1}{2}\Delta\cos^2\frac{\alpha}{2},
\end{eqnarray*}
we obtain
\begin{eqnarray}\label{3.23}
\frac{1}{2}\Delta\log\cos^2\frac{\alpha}{2}&=&\frac{2\mbox{Re}(\overline{a^i_1}\,a^i_{1\bar{1}1})-|a^1_{11}+\overline{a^1_{1\bar{1}}}|^2+|a^i_{11}|^2+|a^i_{1\bar{1}}|^2}{\cos^2\frac{\alpha}{2}}\nonumber\\&{}&+K-\Omega^1_1/\varphi\wedge\bar{\varphi}\nonumber\\
&=&\frac{2\mbox{Re}(\overline{a^i_1}\,a^i_{1\bar{1}1})-2\mbox{Re}(a^1_{11}a^1_{1\bar{1}})+|a^2_{11}|^2+|a^2_{1\bar{1}}|^2}{\cos^2\frac{\alpha}{2}}\nonumber\\&{}&+K-\Omega^1_1/\varphi\wedge\bar{\varphi}.
\end{eqnarray}
Similarly, we also have
\begin{eqnarray}\label{3.24}
\frac{1}{2}\Delta\log\sin^2\frac{\alpha}{2}&=&\frac{2\mbox{Re}(\overline{a^i_{\bar{1}}}\,a^i_{\bar{1}1\bar{1}})-2\mbox{Re}(a^2_{\bar{1}1}\,a^2_{\bar{1}\,\bar{1}})+|a^1_{\bar{1}1}|^2+|a^1_{\bar{1}\,\bar{1}}|^2}{\sin^2\frac{\alpha}{2}}\nonumber\\&{}&+K+\Omega^2_2/\varphi\wedge\bar{\varphi}.
\end{eqnarray}

On the other hand, the Chern minimal condition (\ref{2.27}) and together with (2.18), we have
\begin{equation}\label{3.25}
a^i_{\bar{1}1}=-a^i_{1\bar{1}}=L^i_{jk}\,a^j_1a^k_{\bar{1}},
\end{equation}
which gives
\begin{eqnarray}\label{3.26}
a^i_{\bar{1}1\bar{1}}&=&(L^i_{jk}\,a^j_1a^k_{\bar{1}})_{\bar{1}}\nonumber\\
&=&L^i_{jk\ell}\,a^\ell_{\bar{1}}a^j_1a^k_{\bar{1}}+L^i_{jk\bar{\ell}}\,\overline{a^\ell_1}a^j_1a^k_{\bar{1}}+L^i_{jk}\,a^j_{1\bar{1}}a^k_{\bar{1}}+L^i_{jk}\,a^j_1a^k_{\bar{1}\,\bar{1}}\nonumber\\
&=&L^i_{122}\,\cos\frac{\alpha}{2}\sin^2\frac{\alpha}{2}+L^i_{12\bar{1}}\,\cos^2\frac{\alpha}{2}\sin\frac{\alpha}{2}+L^i_{12}\,a^1_{1\bar{1}}\sin\frac{\alpha}{2}+L^i_{12}\,a^2_{\bar{1}\,\bar{1}}\cos\frac{\alpha}{2},\hspace{0.5cm}
\end{eqnarray}
where $L^i_{jk\ell}$, $L^i_{jk\bar{\ell}}$ defined by
\begin{equation}\label{3.27}
L^i_{jk\ell}\,\omega^\ell+L^i_{jk\bar{\ell}}\,\overline{\omega^\ell}:=dL^i_{jk}-L^i_{\ell k}\,\omega^\ell_j-L^i_{j\ell}\,\omega^\ell_k+L^\ell_{jk}\,\omega_\ell^i.
\end{equation}
Thus, by using (\ref{3.25}) again, we have
\begin{eqnarray}\label{3.28}
2\mbox{Re}(\overline{a^i_{\bar{1}}}\,a^1_{\bar{1}1\bar{1}})&=&\overline{a^i_{\bar{1}}}\,a^i_{\bar{1}1\bar{1}}+a^i_{\bar{1}}\,\overline{a^i_{\bar{1}1\bar{1}}}\nonumber\\
&=&\sin^2\frac{\alpha}{2}\big[(L^2_{122}+\overline{L^2_{122}})\,\cos\frac{\alpha}{2}\sin\frac{\alpha}{2}+(L^2_{12\bar{1}}+\overline{L^2_{12\bar{1}}})\,\cos^2\frac{\alpha}{2}\big]\nonumber\\
&{}&-\cos\frac{\alpha}{2}\sin^3\frac{\alpha}{2}(L^1_{12}L^2_{12}+\overline{L^1_{12}}\,\overline{L^2_{12}})+2\mbox{Re}(a^2_{\bar{1}1}\,a^2_{\bar{1}\,\bar{1}}).
\end{eqnarray}
Similar caculation gives
\begin{eqnarray}\label{3.29}
2\mbox{Re}(\overline{a^i_{1}}\,a^1_{1\bar{1}1})
&=&-\cos^2\frac{\alpha}{2}\big[(L^1_{121}+\overline{L^1_{121}})\,\cos\frac{\alpha}{2}\sin\frac{\alpha}{2}+(L^1_{12\bar{2}}+\overline{L^1_{12\bar{2}}})\,\sin^2\frac{\alpha}{2}\big]\nonumber\\
&{}&-\cos^3\frac{\alpha}{2}\sin\frac{\alpha}{2}(L^1_{12}L^2_{12}+\overline{L^1_{12}}\,\overline{L^2_{12}})+2\mbox{Re}(a^1_{11}\,a^1_{1\bar{1}}).
\end{eqnarray}

\begin{Theorem}\label{Theorem3.1}
Let $f$ be a generic Chern minimal immersion from the compact Riemann surface $\Sigma$ into Hermitian surface $M$. Let $P$ denote the sum of the orders of all complex points, and $Q$ denote the sum of the orders of all anticomplex points. Then
\begin{equation*}
P-Q=-f^*c_1(M)[\Sigma],
\end{equation*}
where $c_1(M)$ is the first Chern form of $M$.
\end{Theorem}

\emph{Proof}. By using the identities (\ref{3.17})-(\ref{3.20}) repeatedly, we have
\begin{eqnarray}\label{3.30}
\cos^2\frac{\alpha}{2}\,(|a^1_{\bar{1}1}|^2+|a^1_{\bar{1}\,\bar{1}}|^2)
&=&\cos^2\frac{\alpha}{2}\sin^2\frac{\alpha}{2}\big(|a^1_{\bar{1}1}|^2+|a^1_{\bar{1}\,\bar{1}}|^2+|a^2_{11}|^2+|a^2_{1\bar{1}}|^2\big),
\end{eqnarray}
and
\begin{eqnarray}\label{3.31}
\sin^2\frac{\alpha}{2}\,(|a^2_{11}|^2+|a^2_{1\bar{1}}|^2)
=\cos^2\frac{\alpha}{2}\sin^2\frac{\alpha}{2}\big(|a^1_{\bar{1}1}|^2+|a^1_{\bar{1}\,\bar{1}}|^2+|a^2_{11}|^2+|a^2_{1\bar{1}}|^2\big).
\end{eqnarray}
Substituting (\ref{3.28})-(\ref{3.31}) into (\ref{3.23}), (\ref{3.24}) respectively, we have
\begin{eqnarray}\label{3.32}
\frac{1}{2}\Delta\log\tan^2\frac{\alpha}{2}\!\!&=&\!\!
(L^2_{122}+\overline{L^2_{122}})\,\cos\frac{\alpha}{2}\sin\frac{\alpha}{2}+(L^2_{12\bar{1}}+\overline{L^2_{12\bar{1}}})\,\cos^2\frac{\alpha}{2}\nonumber\\
&{}&+(L^1_{121}+\overline{L^1_{121}})\,\cos\frac{\alpha}{2}\sin\frac{\alpha}{2}+(L^1_{12\bar{2}}+\overline{L^1_{12\bar{2}}})\,\sin^2\frac{\alpha}{2}\nonumber\\
&{}&
-\sqrt{-1}\mbox{Ric}_M/\varphi\wedge\bar{\varphi},
\end{eqnarray}
where $\mbox{Ric}_M:=\sqrt{-1}(\Omega^1_1+\Omega^2_2)$ is the Ricci form of $M$, $\mbox{Ric}_M/\varphi\wedge\bar{\varphi}$ is the coefficient of $f^*\mbox{Ric}_M$ with respect to $\varphi\wedge\bar{\varphi}$. To simplify the long term involving $L^i_{jk\ell}$, $L^i_{jk\bar{\ell}}$ in (\ref{3.32}), we define
\begin{equation}\label{3.33}
\theta_L:=(L^1_{1i}+L^2_{2i})\,\omega^i-(\overline{L^1_{1i}}+\overline{L^2_{2i}})\,\overline{\omega^i},
\end{equation}
which is globally defined on $M$. Taking the exterior derivatives and applying (\ref{3.27}), we have
\begin{eqnarray}\label{3.34}
d(L^1_{1i}\,\omega^i+L^2_{2i}\,\omega^i)&=&(L^1_{1ij}\,\omega^j+L^1_{1i\,\bar{j}}\,\overline{\omega^j})\wedge\omega^i+L^1_{1i}\Theta^i\nonumber\\
&{}&+(L^2_{2ij}\,\omega^j+L^2_{2i\,\bar{j}}\,\overline{\omega^j})\wedge\omega^i+L^2_{2i}\Theta^i.
\end{eqnarray}
So, from (3.33) and using (\ref{2.8}), (\ref{3.7}), we obtain
\begin{eqnarray}\label{3.35}
f^*d\theta_L&=&(L^1_{12j}\,\omega^j+L^1_{12\,\bar{j}}\,\overline{\omega^j})\wedge\omega^2+L^1_{12}\Theta^2
+(L^2_{21j}\,\omega^j+L^2_{21\,\bar{j}}\,\overline{\omega^j})\wedge\omega^1+L^2_{21}\Theta^1\nonumber\\
&{}&-(\overline{L^1_{12j}}\,\overline{\omega^j}+\overline{L^1_{12\,\bar{j}}}\,\omega^j)\wedge\overline{\omega^2}+\overline{L^1_{12}}\,\overline{\Theta^2}
-(\overline{L^2_{21j}}\,\overline{\omega^j}+\overline{L^2_{21\,\bar{j}}}\,\omega^j)\wedge\overline{\omega^1}+\overline{L^2_{21}}\,\overline{\Theta^1}\nonumber\\
&=&\big[(L^1_{121}+L^2_{122})\cos\frac{\alpha}{2}\sin\frac{\alpha}{2}+L^1_{12\,\bar{2}}\sin^2\frac{\alpha}{2}+L^2_{12\bar{1}}\cos^2\frac{\alpha}{2}\big]
\varphi\wedge\bar{\varphi}\nonumber\\
&{}&+\big[(\overline{L^1_{121}}+\overline{L^2_{122}})\cos\frac{\alpha}{2}\sin\frac{\alpha}{2}+\overline{L^1_{12\,\bar{2}}}\sin^2\frac{\alpha}{2}+\overline{L^2_{12\bar{1}}}\cos^2\frac{\alpha}{2}\big]
\varphi\wedge\bar{\varphi}
\end{eqnarray}
where we have used $L^i_{jk}=-L^i_{kj}$ and $\Theta^i=2L^i_{12}\cos\frac{\alpha}{2}\sin\frac{\alpha}{2}\,\varphi\wedge\bar{\varphi}$. Thus, (\ref{3.32}) and (\ref{3.35}) give
\begin{equation}\label{e.36}
\Delta\log\tan\frac{\alpha}{2}\,\varphi\wedge\bar{\varphi}=-\sqrt{-1}f^*\mbox{Ric}_M+f^*d\theta_L.
\end{equation}
Recall that $dA=\sqrt{-1}\varphi\wedge\bar{\varphi}$, it follows from (3.36) we have
\begin{equation}\label{3.37}
-\frac{1}{2\pi}\Delta\log\tan\frac{\alpha}{2}\,dA=-\frac{1}{2\pi}f^*\mbox{Ric}_M-\frac{\sqrt{-1}}{2\pi}f^*d\theta_L.
\end{equation}
Notice that $c_1(M)=\frac{1}{2\pi}[\mbox{Ric}_M]$, together with the identities (\ref{3.2}) and (\ref{3.6}), the integration of (\ref{3.37}) over $\Sigma$ gives what we want.

\hfill$\Box$

An immediate corollary is

\begin{Corollary}
Let $f$ be a generic Chern minimal immersion from the compact Riemann surface $\Sigma$ into Hermitian surface $M$. Then $P=Q$ if the first Chern class is vanishing.
\end{Corollary}

\emph{Proof}. It follows from $f^*c_1(M)[\Sigma]=0$. \hfill$\Box$

\emph{Remark}. The K$\ddot{a}$hler surface version of this corollary firstly obtained by Al Vitter and also by J.G. Wolfson in \cite{[W-89]}. 

\vspace{0.6cm}

\section{Euler Characteristic of Tangent bundle and Normal Bundle} 

We begin with determining the tangent bundle $T\Sigma$ and normal bundle $T^\perp\Sigma$ of Riemann surface $\Sigma$ immersed in Hermitian surface $(M, J, g)$. It follows from (\ref{3.7}) we have 
\begin{equation}\label{4.1}
\cos\frac{\alpha}{2}\,\omega^1+\sin\frac{\alpha}{2}\,\overline{\omega^2}=\varphi,
\end{equation}
and
\begin{equation}\label{4.2}
\sin\frac{\alpha}{2}\,\overline{\omega^1}-\cos\frac{\alpha}{2}\,\omega^2=0.
\end{equation}
As before, let $\{e_1, e_2\}$ be the dual frame of $\{\omega^1, \omega^2\}$. Then, along $\Sigma$, we obtain that the dual of $\{\cos\frac{\alpha}{2}\,\omega^1+\sin\frac{\alpha}{2}\,\overline{\omega^2},\; \sin\frac{\alpha}{2}\,\overline{\omega^1}-\cos\frac{\alpha}{2}\,\omega^2\}$ is $\{\cos\frac{\alpha}{2}\,e_1+\sin\frac{\alpha}{2}\,\overline{e_2},\; \sin\frac{\alpha}{2}\,\overline{e_1}-\cos\frac{\alpha}{2}\,e_2\}$. So, in terms of orthonormal frame, by reading (\ref{4.1}) and (\ref{4.2}), then the tangent space spanned by
\begin{equation}\label{4.3}
\big\{\frac{1}{\sqrt{2}}\big(\cos\frac{\alpha}{2}\,(e_1+\overline{e_1})+\sin\frac{\alpha}{2}\,(e_2+\overline{e_2})\big),
\frac{\sqrt{-1}}{\sqrt{2}}\big(\cos\frac{\alpha}{2}\,(e_1-\overline{e_1})-\sin\frac{\alpha}{2}\,(e_2-\overline{e_2})\big)\big\},
\end{equation}
and the normal space spanned by
\begin{equation}\label{4.4}
\big\{\frac{1}{\sqrt{2}}\big(\sin\frac{\alpha}{2}\,(e_1+\overline{e_1})-\cos\frac{\alpha}{2}\,(e_2+\overline{e_2})\big),
-\frac{\sqrt{-1}}{\sqrt{2}}\big(\sin\frac{\alpha}{2}\,(e_1-\overline{e_1})+\cos\frac{\alpha}{2}\,(e_2-\overline{e_2})\big)\big\}.
\end{equation}
Here, the basis chosen in (\ref{4.3}) and (\ref{4.4}) are agree with the orientation of $M$. It follows from which and together with (\ref{2.28}), (\ref{2.29}) we obtain the connection 1-form $\rho$ of $T\Sigma$:
\begin{eqnarray}\label{4.5}
\rho\!\!&=&\!\!\big\langle\frac{1}{\sqrt{2}}D\big(\cos\frac{\alpha}{2}\,(e_1+\overline{e_1})+\sin\frac{\alpha}{2}\,(e_2+\overline{e_2})\big),
\frac{\sqrt{-1}}{\sqrt{2}}\big(\cos\frac{\alpha}{2}\,(e_1-\overline{e_1})-\sin\frac{\alpha}{2}\,(e_2-\overline{e_2})\big)\big\rangle\nonumber\\
&=&\sqrt{-1}\big(-\cos^2\frac{\alpha}{2}\,\varphi^1_1+\sin^2\frac{\alpha}{2}\,\varphi^2_2-\cos\frac{\alpha}{2}\sin\frac{\alpha}{2}(\varphi^{\bar{2}}_1-\overline{\varphi^{\bar{2}}_1})\big),
\end{eqnarray}
and the normal connection 1-form $\rho^\perp$ of $T^\perp\Sigma$:
\begin{eqnarray}\label{4.6}
\rho^\perp\!\!&=&\!\!\big\langle\frac{1}{\sqrt{2}}D\big(\sin\frac{\alpha}{2}\,(e_1\!+\!\overline{e_1})\!-\!\cos\frac{\alpha}{2}\,(e_2\!+\!\overline{e_2})\big),
-\frac{\sqrt{-1}}{\sqrt{2}}\big(\sin\frac{\alpha}{2}\,(e_1\!-\!\overline{e_1})\!+\!\cos\frac{\alpha}{2}\,(e_2\!-\!\overline{e_2})\big)\big\rangle\nonumber\\
&=&\sqrt{-1}\big(\sin^2\frac{\alpha}{2}\,\varphi^1_1-\cos^2\frac{\alpha}{2}\,\varphi^2_2-\cos\frac{\alpha}{2}\sin\frac{\alpha}{2}(\varphi^{\bar{2}}_1-\overline{\varphi^{\bar{2}}_1})\big).
\end{eqnarray}

Recalling that the Gausse curvature $K$ given in (\ref{2.3}) and the normal curvature $K^\perp$ determined by
\begin{equation}\label{4.7}
d\rho^\perp=-\sqrt{-1}\,K^\perp\,\varphi\wedge\bar{\varphi}.
\end{equation}
The Euler characteristic of the tangent bundle $T\Sigma$ and normal bundle $T^\perp\Sigma$ are defined by
\begin{equation}\label{4.8}
\chi(T\Sigma):=\frac{1}{2\pi}\int_\Sigma K\,dA,\hspace{0.5cm}\chi(T^\perp\Sigma):=\frac{1}{2\pi}\int_\Sigma K^\perp\,dA,
\end{equation}
respectively.

\vspace{0.3cm}

\begin{Lemma}\label{lemma4.1}
Notations as before, we have
\begin{equation*}
2(|a^1_{\bar{1}\,1}|^2+|a^2_{\bar{1}1}|^2)=L^1_{12}L^2_{12}\sin\alpha=\overline{L^1_{12}}\;\overline{L^2_{12}}\sin\alpha.
\end{equation*}
\end{Lemma}

\emph{Proof}. We first write $2(|a^1_{\bar{1}\,1}|^2+|a^2_{\bar{1}1}|^2)$ as
\begin{equation}\label{4.9}
2(|a^1_{\bar{1}\,1}|^2+|a^2_{\bar{1}1}|^2)=2\cos^2\frac{\alpha}{2}|a^1_{\bar{1}\,1}|^2+2\sin^2\frac{\alpha}{2}|a^1_{\bar{1}\,1}|^2+2\cos^2\frac{\alpha}{2}|a^2_{\bar{1}1}|^2+2\sin^2\frac{\alpha}{2}|a^2_{\bar{1}1}|^2.
\end{equation}
By using (\ref{3.18}) and (\ref{3.25}) repeatedly, the first term in (\ref{4.9}) becomes
\begin{eqnarray}\label{4.10}
2\cos^2\frac{\alpha}{2}|a^1_{\bar{1}\,1}|^2&=& 2\cos^3\frac{\alpha}{2}\sin\frac{\alpha}{2} L^1_{12}\,\overline{a^1_{\bar{1}1}}\nonumber\\
&=&2\cos^2\frac{\alpha}{2}\sin^2\frac{\alpha}{2}L^1_{12}\, a^2_{\bar{1}1}\nonumber\\
&=&\cos^2\frac{\alpha}{2}\sin^2\frac{\alpha}{2}L^1_{12}\, L^2_{12}\sin\alpha.
\end{eqnarray}
Notice that $2\cos^2\frac{\alpha}{2}|a^1_{\bar{1}\,1}|^2$ is real, by taking the conjugate, we also have
\begin{equation}\label{4.11}
2\cos^2\frac{\alpha}{2}|a^1_{\bar{1}\,1}|^2=\cos^2\frac{\alpha}{2}\sin^2\frac{\alpha}{2}\overline{L^1_{12}}\, \overline{L^2_{12}}\sin\alpha.
\end{equation}
Do similar calculations, for the remaining terms in (\ref{4.9}), we have
\begin{equation}\label{4.12}
2\sin^2\frac{\alpha}{2}|a^1_{\bar{1}\,1}|^2=\sin^4\frac{\alpha}{2}L^1_{12}\, L^2_{12}\sin\alpha=\sin^4\frac{\alpha}{2}\overline{L^1_{12}}\, \overline{L^2_{12}}\sin\alpha,
\end{equation}
\begin{equation}\label{4.13}
2\cos^2\frac{\alpha}{2}|a^2_{\bar{1}\,1}|^2=\cos^4\frac{\alpha}{2}L^1_{12}\, L^2_{12}\sin\alpha=\cos^4\frac{\alpha}{2}\overline{L^1_{12}}\, \overline{L^2_{12}}\sin\alpha,
\end{equation}
\begin{equation}\label{4.14}
2\sin^2\frac{\alpha}{2}|a^2_{\bar{1}\,1}|^2=\cos^2\frac{\alpha}{2}\sin^2\frac{\alpha}{2}L^1_{12}\, L^2_{12}\sin\alpha=\cos^2\frac{\alpha}{2}\sin^2\frac{\alpha}{2}\overline{L^1_{12}}\, \overline{L^2_{12}}\sin\alpha.
\end{equation}
So, the identities are follow from (\ref{4.9})-(\ref{4.14}).\hfill$\Box$

\vspace{0.3cm}

Now, we can prove

\begin{Theorem}\label{Theorem4.2}
Let $f$ be a generic Chern minimal immersion from the compact Riemann surface $\Sigma$ into Hermitian surface $M$. Let $P$ denote the sum of the orders of all complex points, and $Q$ denote the sum of the orders of all anticomplex points. Then
\begin{equation*}
P+Q=-(\chi(T\Sigma)+\chi(T^\perp\Sigma)),
\end{equation*}
where $\chi(T\Sigma)$, $\chi(T^\perp\Sigma)$ are the Euler characteristic of the tangent bundle, normal bundle respectively.
\end{Theorem}

\emph{Proof}. It follows from (\ref{4.5}) and (\ref{4.6}), we have
\begin{equation}\label{4.15}
\rho^\perp-\rho=\sqrt{-1}(\varphi^1_1-\varphi^2_2).
\end{equation}
Taking the exterior derivative on both sides of (\ref{4.15}), and using the first identity in Proposition 2.1, we have 
\begin{eqnarray}\label{4.16}
d(\rho^\perp-\rho)&=&\sqrt{-1}\,d(\varphi^1_1-\varphi^2_2)\nonumber\\
&=&\sqrt{-1}\,d\,(\omega^1_1-\omega^2_2)+\sqrt{-1}\,d\,\big[(L^1_{1i}-L^2_{2i})\,\omega^i-(\overline{L^1_{1i}}-\overline{L^2_{2i}})\,\overline{\omega^i}\,\big].
\end{eqnarray}
Then, it follows from (\ref{2.3}), (\ref{4.7}) and (\ref{4.16}) we have
\begin{eqnarray}\label{4.17}
(K^\perp-K)\varphi\wedge\bar{\varphi}=-d\,(\omega^1_1-\omega^2_2)-d\,\big[(L^1_{1i}-L^2_{2i})\,\omega^i-(\overline{L^1_{1i}}-\overline{L^2_{2i}})\,\overline{\omega^i}\,\big].
\end{eqnarray}
By using the structure equation (\ref{2.9}) and the identity (\ref{3.8}), we have
\begin{eqnarray}\label{4.18}
-d\,(\omega^1_1-\omega^2_2)&=&2\,\omega^1_2\wedge\omega^2_1+(-\Omega^1_1+\Omega^2_2)\nonumber\\
&=&-2\,\overline{\omega^2_1}\wedge\omega^2_1+(-\Omega^1_1+\Omega^2_2)\nonumber\\
&=&2(|\sin\frac{\alpha}{2}\overline{a^1_{\bar{1}\,\bar{1}}}-\cos\frac{\alpha}{2}a^2_{11}|^2-|\sin\frac{\alpha}{2}\overline{a^1_{\bar{1}\,1}}-\cos\frac{\alpha}{2}a^2_{1\bar{1}}|^2)\varphi\wedge\bar{\varphi}\nonumber\\
&{}&+(-\Omega^1_1+\Omega^2_2)\nonumber\\
&=&2\big(|a^1_{\bar{1}\bar{1}}|^2+|a^2_{11}|^2-|a^1_{\bar{1}1}|^2-|a^2_{1\bar{1}}|^2\big)\varphi\wedge\bar{\varphi}+(-\Omega^1_1+\Omega^2_2),
\end{eqnarray}
where we have used (\ref{3.25}) and (\ref{3.17})-(\ref{3.20}) repeatedly. For the second term in (\ref{4.17}), by the definition of $L^i_{jk\ell}$, $L^i_{jk\bar{\ell}}$ given in (\ref{3.27}) and together with (\ref{2.8}), (\ref{3.7}), we have
\begin{eqnarray}\label{4.19}
-d\big[(L^1_{1i}-L^2_{2i})\;\omega^i\big]&=& -\big(L^1_{1ij}\omega^j\wedge\omega^i+L^1_{1i\bar{j}}\overline{\omega^j}\wedge\omega^i+L^1_{ji}\omega^j_1\wedge\omega^i-L^j_{1i}\omega_j^1\wedge\omega^i+L^1_{1i}\Theta^i\big)\nonumber\\
&{}&+\big(L^2_{2ij}\omega^j\wedge\omega^i+L^2_{2i\bar{j}}\overline{\omega^j}\wedge\omega^i+L^2_{ji}\omega^j_2\wedge\omega^i-L^j_{2i}\omega_j^2\wedge\omega^i+L^2_{2i}\Theta^i\big)\nonumber\\
&=&\big[(L^2_{122}-L^1_{121})\cos\frac{\alpha}{2}\sin\frac{\alpha}{2}-L^1_{12\bar{2}}\sin^2\frac{\alpha}{2}+L^2_{12\bar{1}}\cos^2\frac{\alpha}{2}\big]\varphi\wedge\bar{\varphi}\nonumber\\
&{}&-2L^1_{12}L^2_{12}\sin\alpha\;\varphi\wedge\bar{\varphi}+2\big(L^1_{12}\,\omega^2_1\wedge\omega^1+2L^2_{12}\,\omega^1_2\wedge\omega^2\big).
\end{eqnarray}
Substituting the expression of $\omega^i$, $\omega^2_1$ given in (\ref{3.7}), (\ref{3.14}) respectively into the last term in (\ref{4.19}), and using  (\ref{3.25}), (\ref{3.17})-(\ref{3.20}) repeatedly, we obtain
\begin{eqnarray}\label{4.20}
2\big(L^1_{12}\,\omega^2_1\wedge\omega^1\!+\!2L^2_{12}\,\omega^1_2\wedge\omega^2\big)\!\!&=&\!\!2L^1_{12}\cos\frac{\alpha}{2}(\sin\frac{\alpha}{2}\overline{a^1_{\bar{1}1}}-\cos\frac{\alpha}{2}a^2_{1\bar{1}})\varphi\wedge\bar{\varphi}\nonumber\\
&{}&+2L^2_{12}\sin\frac{\alpha}{2}(\sin\frac{\alpha}{2}a^1_{\bar{1}1}-\cos\frac{\alpha}{2}\overline{a^2_{1\bar{1}}})\varphi\wedge\bar{\varphi}\nonumber\\
&=&2\big(|a^1_{\bar{1}1}|^2+|a^2_{1\bar{1}}|^2-L^1_{12}\cos^2\frac{\alpha}{2}a^2_{1\bar{1}}+L^2_{12}\sin^2\frac{\alpha}{2}a^1_{\bar{1}1}\big)\varphi\wedge\bar{\varphi}\nonumber\\
&=&2\big(|a^1_{\bar{1}1}|^2\!+\!|a^2_{1\bar{1}}|^2\big)\varphi\wedge\bar{\varphi}\!-\!2\big(L^1_{12}a^2_{1\bar{1}}\!-\!L^1_{12}\sin^2\frac{\alpha}{2}a^2_{1\bar{1}}\big)\varphi\wedge\bar{\varphi}\nonumber\\
&{}&+2\big(L^2_{12}a^1_{\bar{1}1}-L^2_{12}\sin^2\frac{\alpha}{2}a^1_{\bar{1}1}\big)\varphi\wedge\bar{\varphi}\nonumber\\
&=&2L^1_{12}L^2_{12}\sin\alpha\,\varphi\wedge\bar{\varphi}.
\end{eqnarray}
So, it follows from (\ref{4.19}) and (\ref{4.20}), we have
\begin{eqnarray*}
-d\big[(L^1_{1i}-L^2_{2i})\;\omega^i\big]\!\!=\!\!\big[(L^2_{122}-L^1_{121})\cos\frac{\alpha}{2}\sin\frac{\alpha}{2}-L^1_{12\bar{2}}\sin^2\frac{\alpha}{2}+L^2_{12\bar{1}}\cos^2\frac{\alpha}{2}\big]\varphi\wedge\bar{\varphi},
\end{eqnarray*}
and hence
\begin{eqnarray}\label{4.21}
-d\,\big[(L^1_{1i}\!-\!L^2_{2i})\,\omega^i\!-\!(\overline{L^1_{1i}}\!-\!\overline{L^2_{2i}})\,\overline{\omega^i}\big]\!\!\!&=&\!\!\!\big[(L^2_{122}-L^1_{121}+\overline{L^2_{122}}-\overline{L^1_{121}})\cos\frac{\alpha}{2}\sin\frac{\alpha}{2}\nonumber\\
&{}&
-(L^1_{12\bar{2}}\!+\!\overline{L^1_{12\bar{2}}})\sin^2\frac{\alpha}{2}\!+\!(L^2_{12\bar{1}}\!+\!\overline{L^2_{12\bar{1}}})\cos^2\frac{\alpha}{2}\big]\,\varphi\wedge\bar{\varphi},\nonumber\\
&{}&
\end{eqnarray}
Now we substitute (\ref{4.18}) and (\ref{4.21}) back to (\ref{4.17}), we obtain
\begin{eqnarray}\label{4.22}
K^\perp&=&K+2\big(|a^1_{\bar{1}\,\bar{1}}|^2+|a^2_{11}|^2-|a^1_{\bar{1}1}|^2-|a^2_{1\bar{1}}|^2\big)+\big(-\Omega^1_1+\Omega^2_2\big)/\varphi\wedge\bar{\varphi}\nonumber\\
&{}&+\big(L^2_{122}-L^1_{121}+\overline{L^2_{122}}-\overline{L^1_{121}}\big)\cos\frac{\alpha}{2}\sin\frac{\alpha}{2}
-\big(L^1_{12\bar{2}}+\overline{L^1_{12\bar{2}}}\big)\sin^2\frac{\alpha}{2}\nonumber\\
&{}&+\big(L^2_{12\bar{1}}+\overline{L^2_{12\bar{1}}}\big)\cos^2\frac{\alpha}{2}.
\end{eqnarray}

On the other hand, it follows from (\ref{3.23}), (\ref{3.24}), (\ref{3.28})-(\ref{3.31}) we get
\begin{eqnarray}\label{4.23}
\Delta\log\sin\alpha&=&\frac{1}{2}\Delta\log\cos^2\frac{\alpha}{2}+\frac{1}{2}\Delta\log\sin^2\frac{\alpha}{2}\nonumber\\
&=&2K+2\big(|a^1_{\bar{1}\,\bar{1}}|^2+|a^2_{11}|^2+|a^1_{\bar{1}1}|^2+|a^2_{1\bar{1}}|^2\big)+\big(-\Omega^1_1+\Omega^2_2\big)/\varphi\wedge\bar{\varphi}\nonumber\\
&{}&+\big(L^2_{122}-L^1_{121}+\overline{L^2_{122}}-\overline{L^1_{121}}\big)\cos\frac{\alpha}{2}\sin\frac{\alpha}{2}
-\big(L^1_{12\bar{2}}+\overline{L^1_{12\bar{2}}}\big)\sin^2\frac{\alpha}{2}\nonumber\\
&{}&+\big(L^2_{12\bar{1}}+\overline{L^2_{12\bar{1}}}\big)\cos^2\frac{\alpha}{2}-\big(L^1_{12}L^2_{12}+\overline{L^1_{12}}\,\overline{L^2_{12}}\big)\sin\alpha.
\end{eqnarray}
Thus, the identities (\ref{4.22}) and (\ref{4.23}) give
\begin{eqnarray}\label{4.24}
\Delta\log\sin\alpha&=&K+K^\perp+4\big(|a^1_{\bar{1}1}|^2+|a^2_{1\bar{1}}|^2\big)-\big(L^1_{12}L^2_{12}+\overline{L^1_{12}}\,\overline{L^2_{12}}\big)\sin\alpha\nonumber\\
&=&K+K^\perp,
\end{eqnarray}
The last identity holds follows from the Lemma 4.1.
By using (\ref{3.2}) and (\ref{3.6}), the statement follows from the integration of (\ref{4.24}) over Riemann surface $\Sigma$.

\hfill$\Box$

As a corollay, we have
\begin{Corollary}\label{Corollary4.3}
Let $f$ be a generic Chern minimal immersion from the compact Riemann surface $\Sigma$ into Hermitian surface $M$. Let $P$ denote the sum of the orders of all complex points, and $Q$ denote the sum of the orders of all anticomplex points. Then
\begin{eqnarray*}
P=-\big(\chi(T\Sigma)+\chi(T^\perp\Sigma)+f^*c_1(M)[\Sigma]\big)/2,\\
Q=-\big(\chi(T\Sigma)+\chi(T^\perp\Sigma)-f^*c_1(M)[\Sigma]\big)/2.
\end{eqnarray*}
\end{Corollary}

\emph{Proof}. It follows directly from the Theorem \ref{Theorem3.1} and Theorem \ref{Theorem4.2}.

\hfill$\Box$

\vspace{0.2cm}

\section{Applications}

We first use the sum of Gaussian cuvature and normal curvature to characterize the generic real Chern minimal immersion with constant K$\ddot{a}$hler angle. That is

\begin{Theorem}
Let $f$ be a generic Chern minimal immersion from the compact Riemann surface $\Sigma$ into Hermitian surface $M$. Then $f$ is real with constant K$\ddot{a}$hler angle if and only if $K+K^\perp=0$.
\end{Theorem}

\emph{Proof}. Suppose that the K$\ddot{a}$hler angle associated to $f$ is a constant and it belongs to $(0, \pi)$, then the identity (\ref{4.24}) implies $K+K^\perp=0$. Conversely, according to  Theorem 4.2, it is clear that $P+Q= 0$ if $K+K^\perp=0$. This means that $f$ has no complex point and anticomplex point, and hence  $\alpha$ is smooth on $\Sigma$. Applying the maximum principle or divergence theorem to $\Delta\log\sin\alpha= 0$ to get $\sin\alpha$ is a constant.

\hfill$\Box$

The next application involving the topological information of Chern minimal immersion. That is 
\begin{Theorem}
Let $f$ be a generic Chern minimal immersion from the compact Riemann surface $\Sigma$ into compact Hermitian surface $M$. Then, we have
\begin{equation*}
(2-2g)+\big|f^*c_1(M)[\Sigma]\big|+I_f-2D_f\leq -2\min\{P,Q\}\leq 0,
\end{equation*}
where $g$ is the genus of $\Sigma$, $I_f:=\langle D_M^{-1}f_*[\Sigma]\smile D_M^{-1}f_*[\Sigma], [M]\rangle$ is the intersection number of $\Sigma$, $D_f$ is the self-intersection number which  has only regular self intersections of multiplicity $2$.
\end{Theorem}

\emph{Proof}. Let $X$ be the subset of $f(\Sigma)$ in which each point has two preimages. $X$ is a finite subset of $M$. Let $y\in X$ with the preimages $x_1, x_2\in\Sigma$. We evaluate $y$ to have the weight $s(y)=+1$ if $df(T_{x_1}\Sigma)$ and $df(T_{x_2}\Sigma)$ together define the positive orientation on $T_yM$, otherwise $s(y)=-1$. Then, the self-intersection number is defined as 
\begin{equation}\label{5.1}
D_f=\sum\limits_{y\in X}s(y).
\end{equation}
Set $\tilde{X}:=f^{-1}(X)$ and $s(x):=s(f(x))$ for $x\in\tilde{X}$. We define the zero-cycle $[\tilde{X}]\in H_0(\Sigma)$ as follows
\begin{equation*}
[\tilde{X}]:=\sum\limits_{x\in\tilde{X}} s(x)\;x.
\end{equation*}
Then, we have 
\begin{equation}\label{5.2}
[\tilde{X}]=2D_f\cdot a,
\end{equation}
 where $a$ is the generator of $H_0(M)$ dual to $1\in H^0(M)$. Let $D_M:H^*(M)\longrightarrow H_*(M)$ and $D_\Sigma:H^*(\Sigma)\longrightarrow H_*(\Sigma)$ be the Poincare dulity maps, and let $e\in H^2(M)$ be the Euler class of normal bundle $T^\perp\Sigma$. By the result of K.Lashof and S.Smale (\cite{[LS-59]}), R.J.Herbert (\cite{[H-81]}), we also have
\begin{equation}\label{5.3}
[\tilde{X}]=D_\Sigma\big(f^*D_M^{-1}f_*[\Sigma]-e\big),
\end{equation}
where $[\Sigma]\in H_2(\Sigma)$ is the fundamental class of $\Sigma$. Denote the pairing $H^*\times H_*\longrightarrow \textbf{Z}$ by $\langle\;,\;\rangle$, then we have
\begin{eqnarray}\label{5.4}
\langle 1, D_\Sigma f^*D_M^{-1}f_*[\Sigma]\rangle &=&\langle f^*D_M^{-1}f_*[\Sigma], [\Sigma]\rangle\nonumber\\
&=&\langle D_M^{-1}f_*[\Sigma], f_*[\Sigma]\rangle\nonumber\\
&=&\langle D_M^{-1}f_*[\Sigma]\smile D_M^{-1}f_*[\Sigma], [M]\rangle\nonumber\\
&=&I_f.
\end{eqnarray}
On the other hand, by using (\ref{5.3}) and the fact $D_\Sigma e=\chi(T^\perp\Sigma)\cdot a$, we have
\begin{eqnarray}\label{5.5}
\langle 1, D_\Sigma f^*D_M^{-1}f_*[\Sigma]\rangle &=&\langle 1, [\tilde{X}]+D_\Sigma e\rangle\nonumber\\
&=&2D_f+\chi(T^\perp\Sigma).
\end{eqnarray}
Thus, it follows from (\ref{5.4}) and (\ref{5.5}), we obtain
\begin{equation}\label{5.6}
\chi(T^\perp\Sigma)=I_f-2D_f.
\end{equation}
Notice that $\chi(T\Sigma)=2-2g$, so the statement follows from (\ref{5.6}) and the Thoerem 3.1 and Theorem 4.2.

\hfill$\Box$

\begin{Corollary}\label{Corollary5.3}
Let $f$ be a generic Chern minimal embedding from the compact Riemann surface $\Sigma$ into compact Hermitian surface $M$. Then, we have
\begin{equation*}
(2-2g)+\big|f^*c_1(M)[\Sigma]\big|+I_f\leq -2\min\{P,Q\}\leq 0.
\end{equation*}
\end{Corollary}

\emph{Proof}. It follows from $D_f=0$. \hfill$\Box$

\begin{Corollary}\label{Corollary5.4}
Let $M$ be a compact Hermitian surface, $a\in H_2(M,\textbf{Z})$ satisfying $\langle D_M^{-1}a\smile D_M^{-1}a, [M]\rangle\geq 2g_0-1$. Then, $a$ can not be represented by an embedded Chern minimal surface of genus $g\leq g_0$.
\end{Corollary}

\emph{Proof}. If not, suppose $a$ can be represented by an embedded Chern minimal surface of genus $g\leq g_0$, then $I_f\geq 2g_0-1$. This means that 
$(2-2g)+\big|f^*c_1(M)[\Sigma]\big|+I_f\geq 1$, which is a contradiction with the Corollary 5.3.

 \hfill$\Box$

Once the Webester type formulae are obtained, the proof of the Theorem 5.2 and its corollaries are the same as the one of that $M$ is $CP^2$ given in \cite{[EGT-85]} or $M$ is a K$\ddot{a}$hler surface given in \cite{[W-89]}. We write down the details for completeness.

\hspace{0.6cm}

\noindent \textbf{Acknowledgments}.This project is supported by the NSFC (No.11871445), the project of Stable Support for Youth Team in Basic Research Field, CAS(YSBR-001) and the Fundamental Research Funds for the Central Universities. 


\end{document}